# SPACELIKE HELICES IN MINKOWSKI 4-SPACE $E_1^4$

**Mehmet ÖNDER, Hüseyin KOCAYİĞİT, Mustafa KAZAZ**
*Celal Bayar University, Faculty of Science and Art, Department of Mathematics, 45047 Manisa, Turkey.* e-mail: mehmet.onder@bayar.edu.tr

**Abstract**
In this paper, we give some characterizations for spacelike helices in Minkowski space-time $E_1^4$. We find the differential equations characterizing the spacelike helices and also give the integral characterizations for these curves in Minkowski space-time $E_1^4$.

## 1. Introduction

Helix is one of the most fascinating curves in Science and Nature. These curves can be seen in many subjects of Science such as nanosprings, carbon nanotubes, $\alpha$-helices, DNA double and collagen triple helix, lipid bilayers, bacterial flagella in Escherichia coli and Salmonella, aerial hyphae in actinomycetes, bacterial shape in spirochetes, horns, tendrils, vines, screws, springs, helical staircases and sea shells (helico-spiral structures) (see [6,9,14]). Furthermore, in computer aided design and computer graphics, helices can be used for the tool path description, the simulation of kinematic motion or the design of highways, etc. [15]. These curves also appear as solutions of a variational principle associated with an energy density that linearly depend on both curvature and torsion[1,3].

A general helix in Euclidean 3-space $E^3$ is defined by the property that the tangent makes a constant angle with a fixed straight line (the axis of the general helix)[2]. Therefore, a general helix can be equivalently defined as one whose tangent indicatrix is a planar curve. Certainly, the helices in $E^n$ correspond with those whose unit tangent indicatrices are contained in hyperplanes. A classical result for the helices stated by M. A. Lancret in 1802 and first proved by B. de Saint Venant in 1845 (see [12] for details) is: *A necessary and sufficient condition that a curve to be a general helix is that the ratio of the first curvature to the second curvature be constant i.e., $k_1/k_2$ is constant along the curve, where $k_1$ and $k_2$ denote the first and second curvatures of the curve, respectively.* Similar characterization of the helices in Euclidean 4-space $E^4$ is given by Magden: *A curve $x(s)$ is a helix in Euclidean 4-space $E^4$ if and only if the function*

$$\frac{k_1^2}{k_2^2}+\left[\frac{1}{k_3}\frac{d}{ds}\left(\frac{k_1}{k_2}\right)\right]^2 \tag{1}$$

*is constant, where $k_1, k_2$ and $k_3$ are first, second and third curvatures of Euclidean curve $x(s)$, respectively*[11]. Clearly, (1) has a meaning only if $k_1, k_2$ and $k_3$ are nowhere zero, and it is only under this precondition that (1) is a necessary and sufficient condition for a curve of constant slope or helix. Furthermore, some mathematicians studied on helices in Minkowski space. Ilarslan and Boyacioglu defined the position vectors of a timelike and a null helix in Minkowski 3-space $E^3$[8]. In [10], Kocayigit and Onder have given the characterizations for timelike helices in Minkowski 4-space $E_1^4$ [6].

General helices in Lorentzian spaces are curves whose tangent indicatrices are plane curves. It allows one to distinguish between degenerate and non-degenerate helices. In particular, (see [4,5]) helices in $E_1^3$ were geometrically characterized as geodesics of at surfaces. Namely, geodesics in either a right cylinder if the helix is non-degenerate or a $B$-scroll if the helix is degenerate, (see [7] for details on $B$-scrolls).

In this paper, we give some characterizations for spacelike curves to be a helix in Minkowski 4-space $E_1^4$. We show that, differently from the Euclidean case and the case for the timelike helix in $E_1^4$, there are two differential equations characterizing the spacelike helices in $E_1^4$. Furthermore, we give the integral characterizations for spacelike helices in Minkowski 4-space $E_1^4$.

**2. Preliminaries**

Minkowski space-time $E_1^4$ is an Euclidean space $E^4$ provided with the standard flat metric given by
$$g = -dx_1^2 + dx_2^2 + dx_3^2 + dx_4^2$$
where $(x_1, x_2, x_3, x_4)$ is a rectangular coordinate system in $E_1^4$.

Since $g$ is an indefinite metric, recall that a vector $v \in E_1^4$ can have one of three causal characters; it can be spacelike if $g(v,v) > 0$ or $v = 0$, timelike if $g(v,v) < 0$ and null(lightlike) if $g(v,v) = 0$ and $v \neq 0$. Similarly, an arbitrary curve $x(s)$ in $E_1^4$ can locally be spacelike, timelike or null (lightlike), if all of its velocity vectors $x'(s)$ are respectively spacelike, timelike or null (lightlike). Also recall that the pseudo-norm of an arbitrary vector $v \in E_1^4$ is given by $\|v\| = \sqrt{|g(v,v)|}$. Therefore $v$ is a unit vector if $g(v,v) = \pm 1$. The velocity of the curve $x(s)$ is given by $\|x'(s)\|$. Next, the vectors $v, w$ in $E_1^4$ are said to be orthogonal if $g(v,w) = 0$. We say that a timelike vector is *future pointing* or *past pointing* if the first compound of the vector is positive or negative, respectively.

Denote by $\{T(s), N(s), B_1(s), B_2(s)\}$ the moving Frenet frame along the curve $x(s)$ in the space $E_1^4$. Then $T, N, B_1, B_2$ are the tangent, the principal normal, the first binormal and the second binormal fields, respectively. A timelike(spacelike) curve $x(s)$ is said to be parameterized by a pseudo-arclength parameter $s$, i.e. $g(x'(s), x'(s)) = -1$ ($g(x'(s), x'(s)) = 1$).

Let $x(s)$ be a spacelike curve in Minkowski space-time $E_1^4$, parameterized by arclength function of $s$. Then for the curve $x(s)$ the following Frenet equation is given as follows
$$\begin{bmatrix} T' \\ N' \\ B_1' \\ B_2' \end{bmatrix} = \begin{bmatrix} 0 & k_1 & 0 & 0 \\ -\varepsilon_1 k_1 & 0 & k_2 & 0 \\ 0 & \varepsilon_2 k_2 & 0 & k_3 \\ 0 & 0 & \varepsilon_1 k_3 & 0 \end{bmatrix} \begin{bmatrix} T \\ N \\ B_1 \\ B_2 \end{bmatrix}$$
where
$$g(T,T) = 1, \ g(N,N) = \varepsilon_1, \ g(B_2, B_2) = \varepsilon_2, \ g(B_1, B_1) = -\varepsilon_1 \varepsilon_2, \ \varepsilon_1 = \pm 1, \ \varepsilon_2 = \pm 1$$
and recall that the functions $k_1 = k_1(s)$, $k_2 = k_2(s)$ and $k_3 = k_3(s)$ are called the first, the second and the third curvature of the spacelike curve $x(s)$, respectively. We will assume throughout this work that all the three curvatures satisfy $k_i(s) \neq 0$, $1 \leq i \leq 3$. Here, the signs of $\varepsilon_1$ and $\varepsilon_2$ are changed by a rule. The signature rule between $\varepsilon_1$ and $\varepsilon_2$ can be given as follows

| if $\varepsilon_1$ | | then $\varepsilon_2$ | |
|---|---|---|---|
| | +1 | | +1 or -1 |
| | -1 | | +1 |

or

| if $\varepsilon_2$ | | then $\varepsilon_1$ | |
|---|---|---|---|
| | +1 | | +1 or -1 |
| | -1 | | +1 |

(See [12] for details).

### 3. Characterizations of Spacelike Helices in Minkowski 4-space $E_1^4$

Let $x(s)$ be a spacelike curve in Minkowski space-time $E_1^4$, parameterized by arclength function $s$. Similar to the Euclidean case, a spacelike helix can be defined according to the tangent indicatrix: A spacelike curve in $E_1^4$ is said to be a general helix if its unit tangent indicatrix is contained in a hyperplane. Now, a general helix will be called degenerate (a helix with a null axis) or non-degenerate (a helix with a non-null axis) according to whether the corresponding hyperplane is degenerate or non-degenerate, respectively[4,5]. Geometrically, this class of curves can be described as the geodesics in right cylinders with cross sections being curves in hyperplanes.

**Theorem 3.1.** *A regular spacelike curve $x(s)$ is a spacelike helix in Minkowski 4-space $E_1^4$ if and only if the function*

$$\left(\frac{k_1}{k_2}\right)^2 - \varepsilon_1 \frac{1}{k_3^2}\left(\left(\frac{k_1}{k_2}\right)'\right)^2$$

*is constant.*

**Proof.** Let $x(s)$ be a spacelike helix in $E_1^4$. Then $x(s)$ is degenerate(non-degenerate) helix if the axis of the helix is a null(non-null) vector. Let the axis of the spacelike helix $x(s)$ be the unit vector $U$. Then we have

$$g(T,U) = constant.$$

Differentiating this equation with respect to $s$ and using the Frenet formulae we get

$$g(N,U) = 0.$$

Therefore $U$ is in the subspace $Sp\{T, B_1, B_2\}$ and can be written as follows

$$U = \alpha(s)T(s) + \beta(s)B_1(s) + \gamma(s)B_2(s), \qquad (2)$$

where

$$\alpha = -g(U,T) = constant, \quad -\varepsilon_1\varepsilon_2\beta = g(U, B_1), \quad \varepsilon_2\gamma = g(U, B_2).$$

Since $U$ is unit, we have

$$-\alpha^2 - \varepsilon_1\varepsilon_2\beta^2 + \varepsilon_2\gamma^2 = M. \qquad (3)$$

Here $M$ is +1, -1 or 0 depending if $x(s)$ is non-degenerate or degenerate. The differentiation of (2) gives

$$(\alpha k_1 + \varepsilon_2\beta k_2)N + \left(\frac{d\beta}{ds} + \varepsilon_1\gamma k_3\right)B_1 + \left(\frac{d\gamma}{ds} + \beta k_3\right)B_2 = 0,$$

and this equation yields

$$\beta = -\varepsilon_2 \frac{k_1}{k_2}\alpha = -\frac{1}{k_3}\frac{d\gamma}{ds}, \quad \frac{d\beta}{ds} = -\varepsilon_1 \gamma k_3. \tag{4}$$

Since

$$\frac{d\beta}{ds} = -\varepsilon_1 \gamma k_3 \quad \text{and} \quad \frac{d\beta}{ds} = \frac{k_3'}{k_3^2}\frac{d\gamma}{ds} - \frac{1}{k_3}\frac{d^2\gamma}{ds^2},$$

we find the second order linear differential equation in $\gamma$ given by

$$\frac{d^2\gamma}{ds^2} - \frac{k_3'}{k_3}\frac{d\gamma}{ds} - \varepsilon_1 \gamma k_3^2 = 0. \tag{5}$$

In the above equation, if we change variables as $\frac{1}{k_3}\frac{d}{ds} = \frac{d}{dt}$, that is $t = \int_0^s k_3(s)ds$, then we get

$$\frac{d^2\gamma}{dt^2} - \varepsilon_1 \gamma = 0.$$

The last equation has two solutions according to the change of $\varepsilon_1$. When $\varepsilon_1 = +1$, the solution of the equation is

$$\gamma = A\cosh\int_0^s k_3(s)ds + B\sinh\int_0^s k_3(s)ds, \tag{6}$$

and when $\varepsilon_1 = -1$, the solution of the equation is

$$\gamma = A\cos\int_0^s k_3(s)ds + B\sin\int_0^s k_3(s)ds, \tag{7}$$

where $A$ and $B$ are constant.

Let now consider the case $\varepsilon_1 = +1$. Then, from (4) and (6) we have

$$-\left(A\sinh\int_0^s k_3(s)ds + B\cosh\int_0^s k_3(s)ds\right) = -\varepsilon_2 \frac{k_1}{k_2}\alpha = \beta,$$

$$A\cosh\int_0^s k_3(s)ds + B\sinh\int_0^s k_3(s)ds = -\varepsilon_2 \frac{1}{k_3}\left(\frac{k_1}{k_2}\right)'\alpha = \gamma.$$

From these equations it follows that

$$A = \alpha\left(-\varepsilon_2 \frac{k_1}{k_2}\sinh\int_0^s k_3(s)ds + \varepsilon_2 \frac{1}{k_3}\left(\frac{k_1}{k_2}\right)'\cosh\int_0^s k_3(s)ds\right), \tag{8}$$

$$B = \alpha\left(\varepsilon_2 \frac{k_1}{k_2}\cosh\int_0^s k_3(s)ds - \varepsilon_2 \frac{1}{k_3}\left(\frac{k_1}{k_2}\right)'\sinh\int_0^s k_3(s)ds\right). \tag{9}$$

Hence using (3), (8) and (9) we get

$$B^2 - A^2 = \left[\left(\frac{k_1}{k_2}\right)^2 - \frac{1}{k_3^2}\left(\left(\frac{k_1}{k_2}\right)'\right)^2\right]\alpha^2.$$

Therefore

$$\left(\frac{k_1}{k_2}\right)^2 - \frac{1}{k_3^2}\left(\left(\frac{k_1}{k_2}\right)'\right)^2 = constant := m, \tag{10}$$

From (3), (8), (9) and (10) we have
$$B^2 - A^2 = \alpha^2 m = -\varepsilon_2(\alpha^2 + M).$$

Thus, the sign of the constant $m$ agrees with the sign of $B^2 - A^2$. So, if $\varepsilon_2 = +1$ and $U$ is spacelike or ligtlike then, $m$ is negative. If $U$ is timelike then the sign of $m$ depends on $-\varepsilon_2(\alpha^2 - 1)$. In particular, $m = 0$ if and only if $\alpha^2 = 1$. So we can give the following corollary.

**Corollary 3.1.** *Let $x(s)$ be a unit speed spacelike curve with $\varepsilon_1 = 1$ in Minkowski 4-space $E_1^4$ and $U$ be a unit timelike constant vector. Then $\alpha^2 = g(T,U)^2 = 1$ if and only if there exists a constant $D$ such that*
$$\frac{k_1}{k_2}(s) = D\exp\left(\int_0^s k_3(t)dt\right).$$

When $\varepsilon_1 = -1$, by using (7) with similar calculations as above we get that the spacelike curve $x(s)$ with $\varepsilon_1 = -1$ is a helix if and only if

$$\left(\frac{k_1}{k_2}\right)^2 + \frac{1}{k_3^2}\left(\left(\frac{k_1}{k_2}\right)'\right)^2 = constant. \tag{11}$$

Using (10) and (11), we can characterize the spacelike helix $x(s)$ by the property that

$$\left(\frac{k_1}{k_2}\right)^2 - \varepsilon_1 \frac{1}{k_3^2}\left(\left(\frac{k_1}{k_2}\right)'\right)^2 \tag{12}$$

is constant.

Conversely, if the condition (12) is satisfied for a regular spacelike curve we can always find a constant vector $U$ which satisfies the equality $g(T,U) = constant$. It is sufficient to define the unit vector $U$ by

$$U = \left[T - \varepsilon_2 \frac{k_1}{k_2} B_1 + \varepsilon_1\varepsilon_2 \frac{1}{k_3}\left(\frac{k_1}{k_2}\right)' B_2\right].$$

By taking account of the differentiation of (12), differentiation of $U$ gives that $\frac{dU}{ds} = 0$, this means that $U$ is a constant vector. Therefore, the spacelike curve $x(s)$ is a spacelike helix in $E_1^4$.

**Theorem 3.2.** *A regular spacelike curve $x(s)$ in Minkowski space-time $E_1^4$ is a helix if and only if there exists a $C^2$-function $f$ such that*

$$fk_3 = \frac{d}{ds}\left(\frac{k_1}{k_2}\right), \qquad \frac{d}{ds}f(s) = \varepsilon_1 k_3 \frac{k_1}{k_2}. \tag{13}$$

**Proof.** If the spacelike curve $x(s)$ is a helix, from Theorem 3.1 we can write

$$\left(\frac{k_1}{k_2}\right)\frac{d}{ds}\left(\frac{k_1}{k_2}\right) - \varepsilon_1 \frac{1}{k_3}\frac{d}{ds}\left(\frac{k_1}{k_2}\right)\frac{d}{ds}\left[\frac{1}{k_3}\frac{d}{ds}\left(\frac{k_1}{k_2}\right)\right] = 0, \tag{14}$$

and hence we have

$$\frac{1}{k_3}\frac{d}{ds}\left(\frac{k_1}{k_2}\right) = \varepsilon_1 \frac{\left(\frac{k_1}{k_2}\right)\frac{d}{ds}\left(\frac{k_1}{k_2}\right)}{\frac{d}{ds}\left[\frac{1}{k_3}\frac{d}{ds}\left(\frac{k_1}{k_2}\right)\right]}. \tag{15}$$

If we write

$$f(s) = \varepsilon_1 \frac{\left(\frac{k_1}{k_2}\right)\frac{d}{ds}\left(\frac{k_1}{k_2}\right)}{\frac{d}{ds}\left[\frac{1}{k_3}\frac{d}{ds}\left(\frac{k_1}{k_2}\right)\right]}, \tag{16}$$

then we get

$$f(s)k_3 = \frac{d}{ds}\left(\frac{k_1}{k_2}\right). \tag{17}$$

From (15) it can be written

$$\frac{d}{ds}\left[\frac{1}{k_3}\frac{d}{ds}\left(\frac{k_1}{k_2}\right)\right] = \varepsilon_1 k_3 \frac{k_1}{k_2}. \tag{18}$$

By using (17) and (18) we have

$$\frac{d}{ds}f(s) = \varepsilon_1 k_3 \frac{k_1}{k_2}. \tag{19}$$

Conversely, let $f(s)k_3 = \frac{d}{ds}\left(\frac{k_1}{k_2}\right)$. If we define unit vector $U$ given by

$$U = \left[T - \varepsilon_2 \frac{k_1}{k_2} B_1 + \varepsilon_1 \varepsilon_2 f(s) B_2\right] \tag{20}$$

we have $g(U,T) = 1 = const$. Then, $x(s)$ is a spacelike helix with the axis $U$.

**Theorem 3.3.** *A spacelike curve $x(s)$ is a helix in Minkowski space-time $E_1^4$ if and only if the following equality holds,*

$$\frac{k_1}{k_2} = C_1 \eta(\theta) + C_2 \mu(\theta). \tag{21}$$

*Here $C_1$, $C_2$ are constants, $\theta(s) = \int_0^s k_3(s)ds$ and $\eta(\theta) = \cosh(\theta)$, $\mu(\theta) = \sinh(\theta)$, if $\varepsilon_1 = +1$; and $\eta(\theta) = \cos(\theta)$, $\mu(\theta) = \sin(\theta)$, if $\varepsilon_1 = -1$.*

**Proof.** Suppose that spacelike curve $x(s)$ is a helix in Minkowski space-time $E_1^4$. Then the condition in Theorem 3.2 is satisfied. Let us define the $C^2$-function $\theta$ and the $C^1$-functions $m(s)$ and $n(s)$ by

$$\theta(s) = \int_0^s k_3(s)ds, \tag{22}$$

$$\left.\begin{array}{l} m(s) = \dfrac{k_1}{k_2}\eta(\theta) - \varepsilon_1 f(s)\mu(\theta), \\ \\ n(s) = -\varepsilon_1 \dfrac{k_1}{k_2}\mu(\theta) + f(s)\eta(\theta). \end{array}\right\} \tag{23}$$

Here, $\eta(\theta) = \cosh(\theta)$, $\mu(\theta) = \sinh(\theta)$, if $\varepsilon_1 = +1$; and $\eta(\theta) = \cos(\theta)$, $\mu(\theta) = \sin(\theta)$, if $\varepsilon_1 = -1$. If we differentiate equations (23) with respect to $s$ and take account of (22), (17) and (19) we find that $m' = 0$ and $n' = 0$. Therefore, $m(s) = C_1$, $n(s) = C_2$ are constants. Now substituting these in (21) and solving the resulting equations for $\dfrac{k_1}{k_2}$, we get

$$\frac{k_1}{k_2} = C_1\eta(\theta) + C_2\mu(\theta),$$

which is desired equality.

Finally assume that (20) holds. Then from the equations in (23) we get

$$f = C_1\mu(\theta) + C_2\eta(\theta),$$

which satisfies the conditions of Theorem 3.2. So, the spacelike curve $x(s)$ is a spacelike helix in Minkowski space-time $E_1^4$.

**4. Conclusions**

In this paper, the characterizations of spacelike helices are given in Minkowski 4-space $E_1^4$. It is shown that a spacelike curve $x: I \subset IR \to E_1^4$ is a helix if the equation given in (12) holds between the first, second and third curvatures of the curve. Furthermore, the integral characterization of the spacelike helix is given. Geometrically, helix can be described as the geodesics in right cylinders with cross sections being curves in hyperplanes. Thus, the characterizations of spacelike helices can be thought by considering the flat surfaces which contain the helices as geodesics.